\documentclass[12pt]{article}
\usepackage{latexsym}
\usepackage{amsfonts}

\usepackage{amssymb}
\usepackage{epsfig}
\usepackage{amsmath}
\usepackage{amsthm}
\usepackage[mathscr]{eucal}
\usepackage{subfigure}
\usepackage{graphicx}

\newtheorem{Theorem}{Theorem}

\renewcommand{\qed}{\hfill{\ \ \rule{2mm}{2mm}} \vspace{0.2in}}

\begin{document}

\title{Linear Eulerian Extensions of Inhomogenous Random Graphs}
\author{ \textbf{Ghurumuruhan Ganesan}
\thanks{E-Mail: \texttt{gganesan82@gmail.com} } \\
\ \\
IISER Bhopal}
\date{}
\maketitle

\begin{abstract}
The Eulerian extension number of any graph~\(H\) (i.e. the minimum number of edges needed to be added to make~\(H\) Eulerian) is at least~\(t(H),\) half the number of odd degree vertices of~\(H.\) In this paper we consider an inhomogenous random graph~\(G\) whose edge probabilities need not all be the same and use an iterative probabilistic method to obtain sufficient conditions for the Eulerian extension number of~\(G\) to grow \emph{linearly} with~\(t(G).\) We derive our conditions in terms of the average edge probabilities and edge density and also briefly illustrate our result with an example.


\vspace{0.1in} \noindent \textbf{Key words:} Inhomogenous random graphs; linear Eulerian extensions; iterative probabilistic method.

\vspace{0.1in} \noindent \textbf{AMS 2000 Subject Classification:} Primary: 60C05; 60J35;
\end{abstract}

\bigskip

\renewcommand{\theequation}{\arabic{section}.\arabic{equation}}
\setcounter{equation}{0}
\section{Introduction} \label{intro}
For a deterministic graph~\(H,\) Boesch et al. (1977) studied conditions under which~\(H\) is a subeulerian graph; i.e. occurs as a spanning subgraph of an Eulerian graph and also classified all subeulerian graphs. Later Lesniak and Oellermann (1986) presented a detailed survey on subgraphs and supergraphs of Eulerian graphs and multigraphs. For applications of Eulerian extensions, we refer to H\"ohn et al. (2012) and Dorn et al. (2013).

Recently, Ganesan (2020) has studied Eulerian edit numbers of random graphs with uniform edge probabilities, where we show that roughly~\(\frac{n}{4}\) operations suffices to obtain an Eulerian extension with high probability, i.e. with probability converging to one as~\(n \rightarrow \infty.\) The constructed extension is in fact linear, since the number of odd degree vertices is roughly~\(\frac{n}{2}\) with high probability (for details, we refer to Ganesan~(2020)).

In the following section, we consider linear Eulerian extensions of \emph{inhomogenous} random graphs whose edge probabilities need not all be the same. We obtain sufficient conditions for linear extendability in terms of the average vertex degrees and the edge density and use an iterative probabilistic method (described briefly in the Appendix) to construct the desired extension.

\renewcommand{\theequation}{\thesection.\arabic{equation}}
\setcounter{equation}{0}
\section{Linear Eulerian Extensions}\label{sec_lin_eul}
Let~\(K_n\) denote the complete graph with vertex set~\(\{1,2,\ldots,n\}\) and let~\(\{X(h)\}_{h \in K_n}\) be independent random variables with indices over the edge set of~\(K_n\) and satisfying
\begin{equation}\label{x_dist}
\mathbb{P}(X(h) = 1) = p(h)  = 1-\mathbb{P}(X(h) = 0)
\end{equation}
for every edge~\(h.\) If~\(h\) has endvertices~\(u\) and~\(v,\) then we also write~\(p(h) = p(u,v) = p(v,u).\) Let~\(G\) be the random graph formed by the set of all edges~\(h\) satisfying~\(X(h) = 1.\)

Let~\(m(G)  = \#E\) be the number of edges in~\(G.\) For a vertex~\(v\) define~\(d_G(v)\) to be the degree of vertex~\(v,\) i.e., the number of vertices adjacent to~\(v\) in~\(G.\) Let~\(\Delta(G)\) denote the maximum vertex degree in~\(G\) and set~\(\overline{G}\) to be the complement of~\(G,\) i.e.an edge~\((x,y) \in \overline{G}\) if and only if~\((x,y) \notin G.\) A walk in~\(G\) is a sequence of vertices~\({\cal P} := (v_1,v_2,\ldots,v_t)\) such that each~\(v_i\) is adjacent to~\(v_{i+1}.\) If in addition~\(v_l\) is also adjacent to~\(v_1,\) then~\({\cal P}\) is said to be a \emph{circuit}. If each edge~\(G\) occurs exactly once in~\({\cal P},\) then~\({\cal P}\) is an \emph{Eulerian} circuit and~\(G\) is said to be Eulerian. If~\({\cal P}\) has no repeated vertices, then~\({\cal P}\) is said to be a \emph{path}. Similarly if~\({\cal P}\) is a circuit and no vertex in~\({\cal P}\) is repeated, then~\({\cal P}\) is said to be a \emph{cycle}. We say that~\(G\) is \emph{connected} if any two vertices of~\(G\) can be connected by a path and \emph{disconnected} otherwise.

We say that~\(G\) is \emph{Eulerian extendable} or simply extendable if we can convert~\(G\) into an Eulerian graph by simply adding edges from the complement. The resulting graph is then said to be an \emph{Eulerian extension} of~\(G.\) In Ganesan (2020) it is shown that if~\(G_{hom}\) is the homogenous random graph with edge probability~\(p,\) then with high probability~\(G_{hom}\) is extendable and we also obtain estimates on the minimum number of edges needed for such an extension.


If~\(G\) is extendable, then any extension of~\(G\) has only even degree vertices
and so if~\(N_{ext}(G)\) is the minimum number of edges needed to obtain an extension of~\(G,\) then~\(N_{ext}(G) \geq t(G),\) where~\(t(G)\) is \emph{half} the number of odd degree vertices in~\(G.\) Letting~\(E_{disc}\) denote the event that~\(G\) is disconnected, the following result obtains upper bounds for~\(N_{ext}(G)\) in terms the  average edge probabilities
\[\alpha_{low}(n) := \min_{u} \frac{1}{n-1}\sum_{v \neq u }p(u,v),\;\;\;\alpha_{up}(n) := \max_{u} \frac{1}{n-1}\sum_{v \neq u }p(u,v)\] and the average edge density~\(\alpha_e(n) := {{n\choose 2}}^{-1}\sum_{1 \leq u < v \leq n} p(u,v).\)
\begin{Theorem}\label{thm_main} If
\begin{equation}\label{alp_cond2}
\frac{1}{n^{\beta}} \leq \alpha_{low}(n) \leq \alpha_{up}(n) \leq \max\left(\frac{1}{2},1- \sqrt{\frac{\alpha_e(n)}{2}}\right)  - \frac{1}{n^{\gamma}},
\end{equation}
for some constants~\(0 < \beta < \frac{1}{2}\) and~\(0 < \gamma < \frac{1}{2}-\beta,\) then
\begin{equation}\label{n_ext_gam2}
\mathbb{P}\left(N_{ext}(G) \leq 3t(G)\right) \geq 1- e^{-D (\log{n})^2}  - \mathbb{P}(E_{disc}),
\end{equation}
for some constant~\(D > 0.\)
\end{Theorem}
Thus if~\(G\) is connected, then with high probability the number of edges that need to be added to~\(G\) for obtaining an extension, grows linearly in~\(t(G).\)  The condition~\(\alpha_{low}(n) \leq \alpha_e(n) \leq \alpha_{up}(n)\)  is obtained using the fact that the sum of the vertex degrees of a graph is twice the number of edges.  The condition~\(\alpha_e(n) \leq 1-\sqrt{\frac{\alpha_e(n)}{2}} - \frac{1}{n^{\gamma}}\) therefore implies that all values of~\(\alpha_e(n)\) up to~\(\frac{1}{2} - O\left(\frac{1}{n^{\gamma}}\right)\) are permissible.




Before proving Theorem~\ref{thm_main}, we illustrate the condition~(\ref{alp_cond2}) and the event~\(E_{disc}\) with an example. Let~\(0 < b < a < 1\) be two constants, that do not depend on~\(n\) and set~\(p(u,v) =1\) if~\(1 \leq u < v \leq \frac{n}{\log{n}}\) or the edge~\((u,v)\) belongs to the cycle~\((1,2,\ldots,n,1).\) Similarly, set~\(p(u,v) = 0\) if~\(\frac{n}{\log{n}} + 1 \leq u < v \leq \frac{2n}{\log{n}}\) and for any edge~\((n,v)\) containing~\(n\) as an endvertex, we set~\(p(n,v) = a.\) For any other edge~\(e,\) we set~\(p(e) = b.\) For the vertex~\(n,\) we have that
\begin{equation}\label{max_deg_verf}
\sum_{1 \leq v \leq n-1} p(n,v) = a\left(n-1\right)
\end{equation}
and so~\(\alpha_{up}(n) = a.\) Similarly, by considering the vertex~\(1\) we get that
\begin{equation}\label{max_deg_verf2}
\sum_{2 \leq v \leq n} p(1,v) \geq b \left(n-1-\frac{n}{\log{n}}\right)
\end{equation}
and so~\(\alpha_{low}(n) = b(1+o(1)),\)  where~\(o(1) \longrightarrow 0\) as~\(n\rightarrow \infty.\)

The number of edges having an edge probability~\(\geq a\) is at most~\(\left(\frac{n}{\log{n}}\right)^2 + 2n\) and so
\[\sum_{1\leq u<v \leq n} p(u,v) \leq \left(\frac{n}{\log{n}}\right)^2 + 2n + b {n \choose 2}.\] Similarly the number of edges for which the edge probability is zero equals~\(\left(\frac{n}{\log{n}}\right)^2\) and so
\[\sum_{1\leq u<v \leq n} p(u,v)  \geq b \left({n \choose 2} - \left(\frac{n}{\log{n}}\right)^2\right).\]
Combining the above estimates, we get that~\(\alpha_e(n) = b \cdot(1+o(1)).\) Finally, by definition, the cycle~\((1,2,\ldots,n,1)\) on~\(n\) vertices is a subgraph of~\(G\) almost surely and so~\(G\) is connected with probability one.

From the discussions following Theorem~\ref{thm_main}, we see that for every~\((a,b)\) satisfying~\(b< \frac{1}{2}\) and~\(b < a < 1-\sqrt{\frac{b}{2}},\) the random graph~\(G\) has a linear extension in the sense of~(\ref{n_ext_gam2}).





Below, we use the iterative probabilistic method and the following standard deviation estimate, to prove Theorem~\ref{thm_main}. Let~\(X_1,\ldots,X_m\) be independent Bernoulli random variables satisfying~\(\mathbb{P}(X_j = 1) =1-\mathbb{P}(X_j =0) >0.\) If~\(T_m := \sum_{j=1}^{m} X_j\) and~\(\mu_m = \mathbb{E}T_m,\) then for any~\(0 < \epsilon \leq \frac{1}{2}\) we have that
\begin{equation}\label{conc_est_f}
\mathbb{P}\left(|T_m - \mu_m| \geq \epsilon \mu_m\right) \leq \exp\left(-\frac{\epsilon^2}{4} \mu_m\right).
\end{equation}
For a proof of~(\ref{conc_est_f}), we refer to Corollary A.1.14, pp. 312 of Alon and Spencer (2008).

\emph{Proof of Theorem~\ref{thm_main}}: We first consider the case~\(\alpha_{up} \leq \frac{1}{2} - \frac{(\log{n})^3}{n}.\) For vertices~\(u\neq v\) let
\[Y(u,v) := \sum_{\stackrel{1 \leq z \leq n}{z \neq u,v}} (1-X(u,z))(1-X(v,z)) \]
be the number of vertices neither adjacent to~\(u\) nor adjacent to~\(v\) in the graph~\(G.\)
We have that
\begin{eqnarray}
\mathbb{E}Y(u,v)  &=& \sum_{\stackrel{1 \leq z \leq n}{z \neq u,v}} (1-p(u,z))(1-p(v,z)) \nonumber\\
&\geq& n-\sum_{z \neq u,v}p(u,z) - \sum_{z \neq u,v} p(v,z) \nonumber\\
&\geq& n- \sum_{z \neq u}p(u,z) - \sum_{z \neq v} p(v,z) -2 \nonumber\\
&\geq& n(1-2\alpha_{up})-2 \nonumber
\end{eqnarray}

Using~\(\alpha_{up} \leq \frac{1}{2}\left(1-\frac{(\log{n})^3}{n}\right)\) we get that~\(\mathbb{E}Y(u,v) \geq (\log{n})^3-2\) and so applying the standard deviation estimate~(\ref{conc_est_f}) with~\(\epsilon = \frac{1}{\sqrt{\log{n}}},\) we get that
\begin{equation}\label{yuv_est}
\mathbb{P}\left(Y(u,v) \geq \frac{(\log{n})^3}{2}\right) \geq 1-e^{-\frac{(\log{n})^2}{4}}.
\end{equation}
Defining~\(E_{all} := \bigcap_{u \neq v} \left\{Y(u,v) \geq \frac{(\log{n})^3}{2}\right\}\) we get from~(\ref{yuv_est}) and the union bound that
\begin{equation}\label{f_tot_est}
\mathbb{P}(E_{all}) \geq 1-n^2\cdot e^{-\frac{(\log{n})^2}{4}} \geq 1-e^{-\frac{(\log{n})^2}{5} }
\end{equation}
for all~\(n\) large.

Let~\(E_{con}\) denote the event that~\(G\) is connected and suppose that the joint event~\(E_{all} \cap E_{con}\) occurs. Let~\({\cal T} = \{v_1,\ldots,v_{2t}\}, t = t(G),\) be the set of all odd degree vertices of~\(G.\) For each pair~\(\{v_{2i-1},v_{2i}\}, 1 \leq i \leq t,\) pick a vertex~\(q_i\) adjacent to neither~\(v_{2i-1}\) nor~\(v_{2i}\) in~\(G.\) Adding the edges~\(\{(q_i,v_{2i-1}),(q_i,v_{2i})\}_{1 \leq i \leq t}\) to~\(G,\) we then get an Eulerian graph~\(H.\) Each path we have added has exactly two edges and so using~(\ref{f_tot_est}) and the union bound we then get that
\[\mathbb{P}(N_{ext}(G) \geq 2t(G)) \leq \mathbb{P}(E^c_{all} \cup E^c_{con}) \leq e^{-\frac{(\log{n})^2}{5}} + \mathbb{P}(E^c_{con}).\]
Since~\(E^c_{con} = E_{disc},\) this proves~(\ref{n_ext_gam2}) for the case~\(\alpha_e \leq \frac{1}{2} - \frac{1}{n^{\gamma}}\) in~(\ref{alp_cond2}).

We now consider the remaining condition in~(\ref{alp_cond2}) and for notational simplicity, write~\(\alpha_{low}(n),\alpha_{up}(n)\) and~\(\alpha_e(n)\) simply as~\(\alpha_{low},\alpha_{up}\) and~\(\alpha_e,\) respectively. The proof outline is as follows. We first use the conditions in~(\ref{alp_cond2}) to estimate the maximum vertex degree and the number of edges of~\(G.\) Next, we use the iterative probabilistic method to construct the extension of~\(G\) iteratively.

Choose~\(\zeta\) such that~\(\gamma  + \frac{\beta}{2} < \zeta < \frac{1-\beta}{2}\) (this is possible by Theorem statement) and let~\(\epsilon := \frac{1}{n^{\zeta}}.\) Recalling that~\(d(u)\) is the degree of~\(u\) we use the deviation estimate~(\ref{conc_est_f}) and~(\ref{alp_cond2}) to get that
\begin{equation}
\mathbb{P}\left(d(u) \leq \alpha_{up}(1+\epsilon) (n-1)\right) \geq 1- \exp\left(-\frac{\epsilon^2}{4} (n-1)\alpha_{low}\right). \label{temp_ax2}
\end{equation}
Since~\(\alpha_{low} \geq \frac{1}{n^{\beta}}\) (see~(\ref{alp_cond2})), we get that~\(\frac{\epsilon^2}{4} (n-1)\alpha_{low} \geq \frac{1}{8}n^{1-\beta-2\zeta} \geq (\log{n})^2\) for all~\(n\) large, by choice of~\(\zeta.\) Thus the final term in~(\ref{temp_ax2}) is bounded below by~\(1-e^{-(\log{n})^2}\) and by the union bound we get that the maximum vertex degree~\(\Delta(G)\) satisfies
\begin{equation}\label{max_deg_gamma}
\mathbb{P}\left(\Delta(G) \leq \alpha_{up}(1+\epsilon) (n-1)\right) \geq 1-n \cdot \exp\left(-(\log{n})^2\right)
\end{equation}
for all~\(n\) large.

Next we obtain deviation estimates for the number of edges in~\(G.\) Using~(\ref{alp_cond2}) and the lower bound for~\(\alpha_{low}\) in~(\ref{alp_cond2}), the expected number of edges in~\(G\) is
\[\mathbb{E}m(G) = \alpha_e\cdot {n \choose 2} \geq \alpha_{low} \cdot {n \choose 2} \geq n \cdot (\log{n})^2\] for all~\(n\) large and so again using the deviation estimate~(\ref{conc_est_f}) we get that
\begin{equation}\label{max_edge_gamma}
\mathbb{P}\left(m(G) \leq \alpha_e(1+\epsilon) \cdot {n \choose 2}\right) \geq 1- \exp\left(-\frac{\epsilon^2}{4} \cdot n\cdot (\log{n})^2 \right)
\end{equation}
which converges to one as~\(n \rightarrow \infty,\) by choice of~\(\epsilon.\) Defining the event
\begin{equation}\label{e_good_def}
E_{good} := \left\{\Delta(G) \leq \alpha_{up}(1+\epsilon) (n-1)\right\} \bigcap \left\{m(G) \leq \alpha_e(1+\epsilon) \cdot {n \choose 2}\right\},
\end{equation}
we get from~(\ref{max_deg_gamma}),~(\ref{max_edge_gamma}) and the union bound that
\begin{equation}
\mathbb{P}(E_{good}) \geq 1 - n\cdot \exp\left(-(\log{n})^2\right)  - \exp\left(-\frac{\epsilon^2n(\log{n})^2}{4} \right) \label{e_good_est}
\end{equation}
for all~\(n\) large. Letting~\(E_{con}\) denote the event that~\(G\) is connected, we assume henceforth that the event~\(E_{tot} := E_{good}\cap E_{con}\) occurs.

Collect all odd degree vertices of~\(G\) and call it~\({\cal T}.\) Initially all vertices of~\({\cal T}\) are set to be unmarked. Pick any two unmarked vertices~\(x\) and~\(y\) such that the edge~\((x,y) \notin G\) and mark the edge~\((x,y)\) and the endvertices~\(x\) and~\(y.\) Continue this repeatedly until we are left with a set of unmarked vertices~\(\{a_i\}_{1 \leq i \leq 2s},\) that form a clique in~\(G.\)

If there are two vertices~\(a_i, a_j\) with a common neighbour~\(z \notin \bigcup_{i}\{a_i\}\) in the complement graph~\(G,\) then mark the edges~\((a_i,z)\) and~\((a_j,z)\) and also mark~\(a_i\) and~\(a_j.\) Continue this iteratively until we are left with a clique~\({\cal U} :=\{u_i,v_i\}_{1 \leq i \leq t}\) such that any two vertices in~\({\cal U}\) have disjoint neighbourhoods. Add all the marked edges to~\(G\) and call the resultant graph~\(G_{mod}.\) By construction
\begin{equation}\label{del_g_mod}
\Delta(G_{mod}) \leq \Delta(G)+1 \text{ and } m(G_{mod}) \leq m(G) + 2n,
\end{equation}
where we recall that~\(m(G)\) is the number of edges of~\(G.\)

Setting~\(H_0 := G_{mod},\) we now iteratively convert the odd degree vertices in~\(G_{mod}\) into even degree vertices by constructing a increasing sequence of graphs~\(H_0 \subseteq H_1 \subseteq \ldots \subseteq H_t.\) In the graph~\(H_i,\) all the vertices~\(\{u_j,v_j\}_{1 \leq j \leq i}\subseteq {\cal U}\) would have even degrees and so the graph~\(H_t\) obtained at the end would then be the desired extension.

Let~\(\{(Y_i,Z_i)\}_{1 \leq i \leq t}\) be independent and identically distributed (i.i.d.) vertex pairs in~\(\{1,2,\ldots,n\}^2\)  that are also independent of~\(G.\) Let~\(X_i := (Y_i,Z_i)\) and let~\(\mathbb{P}_X\) denote the distribution of~\((X_1,\ldots,X_t).\) Also define~\(E_i\) to be the event that either~\((u_i,Y_i,Z_i,v_i)\) or~\((u_i,Z_i,Y_i,v_i)\) is a path in the complement graph~\(\overline{H}_{i-1}.\) If~\(E_i\) occurs, then we pick exactly one such path and add it to~\(H_{i-1}\) to obtain the new graph~\(H_i.\) Otherwise we simply set~\(H_{i} = H_{i-1}.\)

If~\(F_{tot} := \bigcap_{i=1}^{t}F_i\) occurs then the graph~\(H_t\) obtained at the end of~\(t\) iterations, is the desired Eulerian extension. Our aim is to show that~\(\mathbb{P}_X(F_{tot}) >0\) and we begin by writing~\(F_i = F_{i,1} \cap F_{i,2}\) where~\(F_{i,1}\) be the event that either the edge set~\(\{(u_i,Y_i),(Z_i,v_i)\} \subset \overline{H}_{i-1}\) or~\(\{(u_i,Z_i),(Y_i,v_i)\} \subset \overline{H}_{i-1}\) and~\(F_{i,2}\) is the event that either~\(Y_i=Z_i\) or~\((Y_i,Z_i)\) is an edge of~\(H_{i-1}.\)

Since the vertices in~\(\{u_j,v_j\}_{1 \leq j \leq t}\) form a clique, the degree of the vertices~\(u_{i}\) and~\(v_{i}\) in~\(H_{i-1}\) is the \emph{same} as in~\(G\) and therefore at most~\(n \alpha_{up}(1+\epsilon),\) since~\(E_{good}\) occurs. Thus the~\(\mathbb{P}_X-\)probability that~\(Y_i\) is a neighbour of~\(u_i\) in~\(G\) is at most~\(\alpha_{up}(1+\epsilon).\) Since~\(u_i\) and~\(v_i\) have disjoint neighbourhoods in the complement graph~\(\overline{G},\)  we see that
\begin{equation}\label{pi_est}
\mathbb{P}_X(F_{i,1} \mid X_1,\ldots,X_{i-1})  \geq 2\left(1- \alpha_{up}(1+\epsilon)\right)^2.
\end{equation}
Also, by construction, for any~\(i\) the number of edges in~\(H_{i-1}\) is at most~\[m(G_{mod}) + 3i \leq m(G_{mod}) + 3n \leq m(G) + 4n,\] by~(\ref{del_g_mod}) and so
\begin{equation}\label{qi_est}
\mathbb{P}_X(F_{i,2}^c \mid X_1,\ldots,X_{i-1}) \leq \frac{1}{n} + \frac{2m(G) + 8n}{n^2} \leq \frac{9}{n} + \alpha_e(1+\epsilon).
\end{equation}
Combining~(\ref{pi_est}) and~(\ref{qi_est}), we then get that
\[\mathbb{P}_X\left(F_i\mid X_1,\ldots,X_{i-1}\right) \geq p_i-q_i\] and so proceeding iteratively,
\begin{equation}\label{e_tot_est}
\mathbb{P}_X(F_{tot}) \geq \prod_{i=1}^{t}(p_i-q_i).
\end{equation}

We now use condition~(\ref{alp_cond2}) in the statement of Theorem~\ref{thm_main} to argue that~\(\min_{1 \leq i \leq t}(p_i-q_i) > 0\) strictly. Indeed, since~\(\epsilon = \frac{1}{n^{\zeta}}\) and~\(\frac{1}{n^{\beta}} \leq \alpha_{up}  \leq 1-\sqrt{\frac{\alpha_e}{2}} - \frac{1}{n^{\gamma}},\) we get that
\[1-\alpha_{up}(1+\epsilon) \geq \sqrt{\frac{\alpha_e}{2}}+ \frac{1}{n^{\gamma}} - \frac{1}{n^{\zeta}} \geq \sqrt{\frac{\alpha_e}{2}}+ \frac{1}{2n^{\gamma}}\]  since~\(\gamma  < \zeta\) strictly, by our choice of~\(\zeta.\) Similarly~\(\frac{9}{n} + \alpha_e(1+\epsilon)  \leq \alpha_e + \frac{2}{n^{\zeta}}\) and so from~(\ref{pi_est}) and~(\ref{qi_est}), we get that
\begin{eqnarray}\label{thmulp}
p_i-q_i &\geq& 2\left(\sqrt{\frac{\alpha_e}{2}} + \frac{1}{2n^{\gamma}}\right)^2 -\frac{1}{n}- \alpha_e - \frac{2}{n^{\zeta}} \nonumber\\
&\geq& \frac{\sqrt{2\alpha_e}}{n^{\gamma}} - \frac{1}{n} - \frac{2}{n^{\zeta}} \nonumber\\
&\geq& \frac{\sqrt{2}}{n^{\gamma+\beta/2}} - \frac{1}{n} - \frac{2}{n^{\zeta}},
\end{eqnarray}
by the lower bound in~(\ref{alp_cond2}). Since we have chosen~\(\gamma +\beta/2< \zeta\) strictly,  we get that~\(p_i-q_i>0\) strictly irrespective of~\(i.\) Summarizing, we see that if the event~\(E_{good} \cap E_{con}\) occurs, then~\(N_{ext}(G) \leq 3t(G).\) From the estimate~(\ref{e_good_est}) for the event~\(E_{good}\) and the union bound, we obtain~(\ref{n_ext_gam2}) and this completes the proof of the Theorem.~\(\qed\)








\underline{\emph{Acknowledgements}}: I thank Professors Rahul Roy and C. R. Subramanian for crucial comments and also thank IMSc and IISER Bhopal for my fellowships.

\bibliographystyle{plain}

\begin{thebibliography}{10}
\bibitem{alon} N. Alon and J. Spencer. (2008).
\newblock{\em The Probabilistic Method}.
\newblock{Wiley}.


\bibitem{tind} F. T. Boesch, C. Suffel and R. Tindell. (1977).
\newblock{The Spanning Subgraph of Eulerian Graphs}.
\newblock{Journal of Graph Theory}, \textbf{1}, pp. 79--84.


\bibitem{boll} B. Bollob\'as. (2001).
\newblock {Random Graphs}.
\newblock {Cambridge University Press}.


\bibitem{dorn} F. Dorn, H. Moser, R. Niedermeier and M. Weller. (2013).
\newblock{Efficient Algorithms for extension and Rural Postman Problem}.
\newblock{ SIAM Journal on Discrete Mathematics}, \textbf{27}, pp. 75--94.



\bibitem{gan2} G. Ganesan. (2020).
\newblock{Deviation Estimates for Eulerian Edit Numbers of Random Graphs}.
\newblock{Statistics and Probability Letters}, \textbf{171}, pp. 1--7.

\bibitem{hohn} W. H\"ohn, T. Jacobs, and N. Megow. (2012).
\newblock{On Eulerian extensions and Their Application to No-Wait Flowshop Scheduling}.
\newblock{Journal of Scheduling}, \textbf{15}, pp. 295--309.

\bibitem{lesniak} L. Lesniak and O. R. Oellermann. (1986).
\newblock{An Eulerian Exposition}.
\newblock{Journal of Graph Theory}, \textbf{10} pp. 277--297.


\bibitem{west} D. B. West. (2001).
\newblock{Introduction to Graph Theory}.
\newblock{Prentice Hall}.


\end{thebibliography}

\end{document}